\documentclass[12pt,a4paper]{amsart}

\usepackage[utf8]{inputenc}
\usepackage[T1]{fontenc}
\usepackage[english]{babel}
\usepackage{cite}

\usepackage{indentfirst}
\usepackage{amssymb}
\usepackage{amsfonts}
\usepackage{amsmath}
\usepackage{amsthm}
\usepackage[top=2.5cm, bottom=2.5cm, left=2.5cm, right=2.5cm]{geometry}
\usepackage{amsopn}
\usepackage[colorlinks]{hyperref}

\theoremstyle{plain}
\newtheorem{thm}{Theorem}

\newtheorem{lem}[thm]{Lemma}
\newtheorem{prop}[thm]{Proposition}
\newtheorem{cor}[thm]{Corollary}

\theoremstyle{definition}
\newtheorem{defn}[thm]{Definition}

\newtheorem*{example*}{Example}

\newtheorem*{rem*}{Remark}
\newtheorem{rem}[thm]{Remark}

\newcommand{\ep}{\eta}
\newcommand{\W}{\mathcal{W}}
\newcommand{\R}{\mathbb{R}}

\DeclareMathOperator{\diam}{diam}
\DeclareMathOperator{\supp}{supp}
\DeclareMathOperator{\dist}{dist}

\title{On density of smooth functions in weighted fractional Sobolev spaces}
\author[B{.} Dyda]{Bart{\l}omiej Dyda}
\author[M{.} Kijaczko]{Micha\l{} Kijaczko}


\address[B.D. and M.K.]{Faculty of Pure and Applied Mathematics\\ Wroc{\l}aw University 
	of Science and Technology\\
	Wybrze\.ze Wyspia\'nskiego 27,
	50-370 Wroc{\l}aw, Poland
}
\email{bdyda@pwr.edu.pl}
\email{michal.kijaczko@pwr.edu.pl}
\thanks{B.D. was partially supported by grant NCN 2015/18/E/ST1/00239.}
\thanks{\textcopyright\,2020.\;This manuscript version is made available under the CC-BY-NC-ND 4.0 license \url{http://creativecommons.org/licenses/by-nc-nd/4.0/}}

\begin{document}

\begin{abstract}
We prove that smooth $C^\infty$ functions are dense in weighted fractional Sobolev spaces on an arbitrary open set, under some mild conditions on the weight. We also obtain a~similar result in non-weighted spaces defined by some kernel similar to $x\mapsto |x|^{-d-sp}$. One may consider the results to be a~version of the Meyers--Serrin theorem.
\end{abstract}

	\maketitle
	\tableofcontents 
	
	\section{Introduction}
	We discuss the problem of density of smooth functions in the fractional Sobolev space $W^{s,p}(\Omega)$, as well as in the weighted fractional Sobolev space $W^{s,p}(\Omega, w)$;
	for the definition of the latter we refer the reader to Section~\ref{sec:weighted}.
	It turns out that for weights $w$ which are locally comparable to a constant on $\Omega$ or continuous, and which satisfy certain integrability property \eqref{conditon}, smooth functions $C^\infty(\Omega)$ are dense in $W^{s,p}(\Omega, w)$, see Theorem~\ref{thm:weighted-dense}.
	
	Our strategy of the proof follows the approach of \cite[proof of Theorem 3.25]{McLean}, in that we first decompose the function $f$ being approximated into the sum of functions $f_n$ supported on the (enlarged) Whitney cubes, which is done by using a~partition of unity.
	Then we convolve each $f_n$ with a~dilation of a~fixed smooth function. In the non-weighted case, the scale of the dilation is dependent on the size of the Whitney cube, to make sure that the support of the convolution does not grow too much. That way we obtain a~family of linear operators $P^{\ep_{k}}$, each mapping the function to a~smooth approximating function, with the error of approximation going to zero when $\ep_{k}$ are sufficiently small. In the weighted case, the scale of the dilation is dependent also on the function being approximated, and the resulting approximating operators are no longer linear.
	
	The proof works for general open sets $\Omega\subset \R^d$, and the result seems to be new in the case of the weighted Sobolev spaces, see Theorem~\ref{thm:weighted-dense}, or a more general kernel, see Theorem~\ref{thm:K}.
        The other standard approach to prove such a~density result is to use the extension theorem \cite{MR3280034}, however it does not hold for all open sets $\Omega$.

	Our paper is motivated by the article \cite{MR3420496}, where the authors consider a~similar problem for weights $w(x)=|x|^{-a}$ in $\R^d$ (or, translated to our setting, in $\R^d\setminus\{0\}$).
	They consider however the density of the compactly supported smooth functions, the problem that we do not address. We note here that if one knows that the compactly supported functions (not necessarily smooth) are dense in $W^{s,p}(\Omega, w)$, then our result immediately gives the density of the space $\mathcal{C}_{C}^\infty(\Omega)$, see Proposition~\ref{prop:support}.
	
	Let us also mention other articles on similar topics.
	In \cite{MR3605166} Luiro and Vähäkangas considered slightly different fractional Sobolev spaces, that are equipped with the seminorm 
	$$
	|f|_{W^{s,p,K}(\R^d)}=\left(\int_{\R^d}\int_{\R^d}\frac{|f(x)-f(y)|^{p}}{|x-y|^{sp}}K(x-y)\,dx\,dy\right)^{\frac{1}{p}},
	$$
	where the kernel $K$ does not have to be radial. The authors find some condition which is sufficient for the space $\mathcal{C}^{\infty}(\R^{d}) \cap W^{s,p,K}(\R^{d})$ to be dense in $W^{s,p,K}(\R^{d})$  (see \cite{MR3605166}, (3.8) and Lemma 3.4). We obtain a~similar result, Theorem~\ref{thm:K}, with more general sets $\Omega$, but less general kernels $K$.
	
	In \cite{MR3310082} Fiscella, Servadei and Valdinoci considered similar Sobolev space $X_{0}^{s,p}(\Omega)$ of functions $f$ with the finite norm
	\[
	\|f\|_{L^p(\R^d)} + \left(\int_{\R^d\times\R^d} |f(x)-f(y)|^p K(x-y)\,dx\,dy\right)^{1/p},
	\]
	but vanishing outside $\Omega$, with some assumptions on the kernel $K$. The authors proved that the space $\mathcal{C}_{C}^{\infty}(\Omega)$ of smooth functions that are compactly supported in $\Omega$, is dense in $X_{0}^{s,p}(\Omega)$, when $\Omega$ is either a hypograph or a domain with continuous boundary (see \cite{MR3310082}, Theorems 2 and 6).
	
	In \cite{MR3989177} Baalal and Berghout considered fractional Sobolev spaces with variable exponents $W^{s,q(\cdot),p(\cdot,\cdot)}(\Omega)$ and
	proved that under certain conditions for the functions $p$ and $q$, compactly supported, smooth functions are dense in $W^{s,q(\cdot),p(\cdot,\cdot)}(\Omega)$.
	
	The authors would like to thank Antti V. V\"ah\"akangas and Victor Nistor for helpful discussions on the subject, and the anonymous reviewer for useful comments.
        We have been informed that a~result similar to our Theorem~\ref{thm:K} has been independently obtained by Foghem Gounoue Guy Fabrice, to be published in his Ph.D. thesis.
	
	\section{Operator $P^{\eta}$}
	\subsection{Definition}
	Let $\Omega\subset\R^{d}$ be an open set and $\W=\{Q_{1},Q_{2},\dots\}$ be a Whitney decomposition of $\Omega$ into cubes, like in \cite{MR0290095}. Choose $\varepsilon$ such that $(1+\varepsilon)^{2}<\frac{5}{4}$. Let also $\{\psi_{n}:n\in\mathbb{N}\}$ be a partition of unity, that is $\psi_{n}(x)=1$, when $x\in Q_{n}$, $\psi_{n}=0$ outside $Q_{n}^{*}$, where $Q_{n}^{*}$ is the cube $Q_{n}$ ''blown up'' $1+\varepsilon$ times (the cube with the same center, but the length of the edge $1+\varepsilon$ times longer), $\psi_{n}$ is a class of $\mathcal{C}^{\infty}_{C}$ and $\sum_{n=1}^{\infty}\psi_{n}=1$. Let $p\in[1,\infty)$ and $f\in L^{p}(\Omega)$.
	
	We note that $|\psi_n(x)-\psi_n(y)| \leq \frac{C|x-y|}{l(Q_{n})}\wedge 1$ for some constant $C>0$.
	
	Let us fix a~function $h\colon\R^{d}\to\R$ such that $h\geq 0$, $\int_{\R^{d}}h(x)\,dx=1$, $\supp h = B(0,1)$ and $h\in C^\infty(\R^d)$.
	For $\delta>0$ we define the dilation
	\[
	h_{\delta}(x)=\frac{1}{\delta^{d}}h\left(\frac{x}{\delta}\right),\quad (x\in \R^d).
	\]
	The function $h_{\delta}$ is a class of $\mathcal{C}_{C}^{\infty}(\R^{d})$ and $\int_{\R^{d}}h_{\delta}(x)\,dx=1$ for every $\delta>0.$
	
	For a function $g:\R^d\to \R$ and $t\in \R^d$ we define its translation $\tau_t g$ by the formula
	\[
	\tau_t g(x)=g(x-t), \quad (x\in \R^d).
	\]

	Let $\ep:\W \to (0,\infty)$ be a~function such that $\ep(Q) < \frac{\varepsilon}{2} l(Q)$ for every $Q\in\W$,
	where $l(Q)$ denotes the length of the edge of the cube $Q$. 
	In particular, we may take $\ep=\delta l$, where $\delta\in (0,\varepsilon/2)$.
	For such a function $\ep$ we define the operator $P^{\ep}$ as
	\begin{equation}\label{def:P}
	P^{\ep}f=\displaystyle\sum_{n=1}^{\infty}(f\psi_{n})*h_{\ep(Q_{n})}, \quad f\in L^1_{loc}(\Omega),
	\end{equation}
	where we put $f=0$ on $\R^d\setminus\Omega$.
	
	\begin{prop}\label{prop:cinfty}
		The operator $P^{\ep}$ is well defined and $P^{\ep}f\in\mathcal{C}^{\infty}(\Omega)$ for $f\in L^1_{loc}(\Omega)$.
	\end{prop}
	\begin{proof}
		We observe that the function  $(f\psi_{n})*h_{\ep(Q_{n})}$,
		\begin{displaymath}
		(f\psi_{n})*h_{\ep(Q_{n})}(x)=\int_{\R^{d}}f(x-y)\psi_{n}(x-y)h_{\ep(Q_{n})}(y)\,dy,
		\end{displaymath}
		vanishes outside $Q_n^{**}$.
		Indeed, if $x\notin Q_n^{**}$, then either  $x-y\notin Q_n^*$, which implies $\psi_{n}(x-y)=0$, or $y\notin B(0,\ep(Q_n))$, which implies $h_{\ep(Q_{n})}(y)=0$, because if $x-y\in Q_{n}^{*}$, then $x\in Q_{n}^{*}+y\subset Q_{n}^{**}$ for $|y|<\eta(Q_{n})$, thanks to our choice of $\varepsilon$.
		
		Since $Q_n^{**}\subset \frac{5}{4}Q_n$ by our choice of $\varepsilon$,
		each point $x\in\Omega$ belongs to at most $12^{d}$ cubes $Q_{n}^{**}$ (see \cite{MR0290095}, chapter VI).
		Therefore the sum \eqref{def:P} has at each point only finitely many nonzero terms, thus the result follows.
	\end{proof}

	\begin{prop}\label{prop:support}
	  If $f\in L^1_{loc}(\Omega)$ satisfies $f=0$ outside a~compact set $K\subset \Omega$,
          then also $P^{\ep}(f)=0$ outside some compact set $K'\subset \Omega$.
	\end{prop}
	\begin{proof}
          We observe that only finitely many of the functions $f\psi_n$ are not identically zero.
          Since  $\supp (f\psi_{n})*h_{\ep(Q_{n})} \subset Q_n^{**}$, it follows that $\supp P^{\ep}(f)$
          is contained in a~finite union of cubes $Q_n^{**}$, which is a~compact subset of~$\Omega$.
          \end{proof}
	\subsection{Convergence of the operator $P^{\ep}$ in $L^{p}(\Omega)$}
	\begin{thm}\label{thm:Lp}
		Let $p\in[1,\infty)$ and $f\in L^{p}(\Omega)$. Then
		\begin{displaymath} 
		\displaystyle\lim_{k\to\infty}\|P^{\ep_k}f-f\|_{L^{p}(\Omega)}=0,
		\end{displaymath}
		provided $\lim_{k\to\infty} \ep_k(Q)=0$ for every $Q\in\W$.
	\end{thm}
	
	\begin{proof}
		We have
		\begin{align*}
		\|P^{\ep_k}f-f\|_{L^{p}(\Omega)}^{p}&=\displaystyle\int_{\Omega}|P^{\ep_k}f(x)-f(x)|^{p}dx\\
		&=\displaystyle\int_{\Omega}\left|\sum_{n=1}^{\infty}(f\psi_{n})*h_{\ep_k(Q_{n})}(x)-\sum_{n=1}^{\infty}f(x)\psi_{n}(x)\right|^{p}\,dx\\
		&\leq \displaystyle\int_{\Omega}\left(\sum_{n=1}^{\infty}\left|(f\psi_{n})*h_{\ep_k(Q_{n})}(x)-f(x)\psi_{n}(x)\right|\right)^{p}\,dx.
		\end{align*}
		The sum above is finite at each point $x$ and has at most  $12^{d}$ nonzero terms. Thus, recalling that $\int_{\R^{d}}h_{t}(x)dx=1$ for every $t>0$ and using Jensen inequality we obtain
		\begin{align}
		&\displaystyle\int_{\Omega}\left(\sum_{n=1}^{\infty}\left|(f\psi_{n})*h_{\ep_k(Q_{n})}(x)-f(x)\psi_{n}(x)\right|\right)^{p}\,dx\nonumber\\
		&\leq  M \displaystyle\int_{\Omega}\sum_{n=1}^{\infty}\left|(f\psi_{n})*h_{\ep_k(Q_{n})}(x)-f(x)\psi_{n}(x)\right|^{p}\,dx\nonumber\\
		&=M\displaystyle\int_{\Omega}\sum_{n=1}^{\infty}\left|\int_{\R^{d}}\left(f(x-y)\psi_{n}(x-y)-f(x)\psi_{n}(x)\right)h_{\ep_k(Q_{n})}(y)dy\right|^{p}\,dx\nonumber\\
		&\leq M\displaystyle\sum_{n=1}^{\infty}\int_{\Omega}\int_{\R^{d}}\left|f(x-y)\psi_{n}(x-y)-f(x)\psi_{n}(x)\right|^{p}h_{\ep_k(Q_{n})}(y)\,dy\,dx\nonumber\\
		&=M\sum_{n=1}^{\infty}\int_{\R^{d}}\|\tau_y \left(f\psi_{n}\right)-f\psi_{n}\|^{p}_{L^{p}(\R^{d})}h_{\ep_k(Q_{n})}(y)\,dy\nonumber\\
		&=M\sum_{n=1}^{\infty}\int_{\R^{d}}\|\tau_{\ep_k(Q_{n})u}\left(f\psi_{n}\right)-f\psi_{n}\|^{p}_{L^{p}(\R^{d})}h(u)\,du,\label{fpsin}
		\end{align}
		where $M=12^{d(p-1)}$. Furthermore,
		\[
		\int_{\R^{d}}\|\tau_{\ep_k(Q_{n})u}\left(f\psi_{n}\right)-f\psi_{n}\|^{p}_{L^{p}(\R^{d})}h(u)\,du\leq 2^{p}\|f\psi_{n}\|^{p}_{L^{p}(\R^{d})},
		\]
		and 
		\[
		\sum_{n=1}^{\infty}\|f\psi_{n}\|^{p}_{L^{p}(\R^{d})}=\sum_{n=1}^{\infty}\int_{Q_{n}^{*}}|f(x)\psi_{n}(x)|^{p}\,dx\leq 12^{d}\|f\|^{p}_{L^{p}(\R^{d})}<\infty.
		\]
		Since $\displaystyle\lim_{k\to\infty}\|\tau_{\ep_k(Q_{n})u}\left(f\psi_{n}\right)-f\psi_{n}\|^{p}_{L^{p}(\R^{d})}=0$, using Lebesgue dominated convergence theorem twice in (\ref{fpsin}) we get the assertion of the Theorem.
	\end{proof}
	
	\section{Sobolev spaces}
	For a~measurable function $f$ defined on $\Omega \subset \R^d$, we define its \emph{Gagliardo seminorm} by
	\[
	[f]_{W^{s,p}(\Omega)} = \left( \int_\Omega \int_\Omega \frac{|f(x)-f(y)|^p}{|x-y|^{d+sp}}\,dy\,dx \right)^{1/p}.
	\]
	For $0<s<1$ and $1\leq p<\infty$ we define the \emph{fractional Sobolev space} $W^{s,p}(\Omega)$ as
	\[
	W^{s,p}(\Omega) = \{ f\in L^p(\Omega): [f]_{W^{s,p}(\Omega)} < \infty \}.
	\]
	
	\subsection{Convergence of the  operator $P^{\ep}$ in Gagliardo seminorm}
	
	\begin{lem}\label{lem:gn}
		Suppose that $\Omega\subset\R^{d}$ and $f\in W^{s,p}(\Omega)$.
		Then
		\begin{equation}\label{eq:Pandgn}
		[P^{\ep_k}f-f]^{p}_{W^{s,p}(\Omega)} \leq M\sum_{n=1}^{\infty}\int_{\R^{d}}\|\tau_{\ep_k(Q_{n})u}\left(g_{n}\right)-g_{n}\|^{p}_{L^{p}(\R^{2d})}h(u)\,du,
		\end{equation}
		where $M=12^{d(p-1)}$, and
		\begin{equation}\label{def:gn}
		g_{n}(x,y)= \begin{cases} \dfrac{f(x)\psi_{n}(x)-f(y)\psi_{n}(y)}{|x-y|^{\frac{d}{p}+s}}\,, & \quad x,y\in\Omega;\\
		0, & \quad (x,y) \in (\R^d\times \R^d )\setminus (\Omega \times \Omega).
		\end{cases}
		\end{equation}
		Furthermore,
		\begin{equation}
		\|g_n\|^p_{L^p(\R^{2d})} \leq c(p,d,s) \left( [f]^p_{W^{s,p}(Q_n^*)} + \|f\|^{p}_{L^{p}(Q_{n}^{*})}l(Q_{n})^{-sp} \right) < \infty
		\end{equation}
		for some constant $c(p,d,s)$ depending only on $p$, $d$, $s$.
	\end{lem}
	
	\begin{proof}
		By arguments similar to that from the proof of Theorem \ref{thm:Lp},
		\begin{align}
		[P^{\ep_k}f&-f]^{p}_{W^{s,p}(\Omega)}=\int_{\Omega}\int_{\Omega}\frac{|P^{\ep_k}f(x)-f(x)-P^{\ep_k}f(y)+f(y)|^{p}}{|x-y|^{d+sp}}\,dx\,dy\nonumber\\
		&\leq M\displaystyle\sum_{n=1}^{\infty}\int_{\Omega}\int_{\Omega}\int_{\R^{d}}\frac{\left|(f\psi_{n})(x-t)-(f\psi_{n})(x)-(f\psi_{n})(y-t)+(f\psi_{n})(y)\right|^{p}}{|x-y|^{d+sp}}\nonumber\\
		&\qquad\times h_{\ep_k(Q_{n})}(t)\,dt\,dx\,dy\nonumber\\
		&\leq M\sum_{n=1}^{\infty}\int_{\R^{d}}\|\tau_t\left(g_{n}\right)-g_{n}\|^{p}_{L^{p}(\R^{2d})}h_{\ep_k(Q_{n})}(t)\,dt\label{gn0}\\
		&=M\sum_{n=1}^{\infty}\int_{\R^{d}}\|\tau_{\ep_k(Q_{n})u}\left(g_{n}\right)-g_{n}\|^{p}_{L^{p}(\R^{2d})}h(u)\,du\label{gn},
		\end{align}
		which proves the first part of the Lemma.
		To prove the remaining part, we observe that
		\begin{align*}
		|f(x)\psi_{n}(x)-f(y)\psi_{n}(y)|^{p}&=|f(x)\psi_{n}(x)-f(x)\psi_{n}(y)+f(x)\psi_{n}(y)-f(y)\psi_{n}(y)|^{p}\\
		&\leq 2^{p-1}(|f(x)|^{p}|\psi_{n}(x)-\psi_{n}(y)|^{p}+|\psi_{n}(y)|^{p}|f(x)-f(y)|^{p}).
		\end{align*}
		Since $\supp\psi_{n}\subset Q_{n}^{*}$, 
		\begin{align*}
		\|g_{n}&\|_{L^{p}(\R^{2d})}^{p}=\int_{\Omega}\int_{\Omega}\frac{|f(x)\psi_{n}(x)-f(y)\psi_{n}(y)|^{p}}{|x-y|^{d+sp}}\,dx\,dy\\
		&\leq 2\int_{\Omega}\int_{Q_{n}^{*}}\frac{|f(x)\psi_{n}(x)-f(y)\psi_{n}(y)|^{p}}{|x-y|^{d+sp}}\,dx\,dy\\
		&\leq 2^{p}\int_{\Omega}\int_{Q_{n}^{*}}\frac{|f(x)|^{p}|\psi_{n}(x)-\psi_{n}(y)|^{p}}{|x-y|^{d+sp}}\,dx\,dy+2^{p}\int_{\Omega}\int_{Q_{n}^{*}}\frac{|\psi_{n}(y)|^{p}|f(x)-f(y)|^{p}}{|x-y|^{d+sp}}\,dx\,dy\\
		&=:2^p(I_1+I_2).
		\end{align*}
		We have $|\psi_{n}(y)|\leq 1$, thus 
		\begin{displaymath}
		I_2\leq \int_{Q_{n}^{*}}\int_{Q_{n}^{*}}\frac{|f(x)-f(y)|^{p}}{|x-y|^{d+sp}}\,dx\,dy = [f]^{p}_{W^{s,p}(Q_n^*)}<\infty.
		\end{displaymath}
		Since $|\psi_{n}(x)-\psi_{n}(x+w)|\leq\frac{C|w|}{l(Q_{n})}\wedge 1$, therefore
		\begin{align*}
		I_1&=\int_{Q_{n}^{*}}\int_{\Omega-x}\frac{|f(x)|^{p}|\psi_{n}(x)-\psi_{n}(x+w)|^{p}}{|w|^{d+sp}}\,dw\,dx\\
		&\leq \int_{Q_{n}^{*}}|f(x)|^{p}\int_{\Omega-x}\left(\frac{C^{p}|w|^{p}}{l(Q_{n})^{p}}\wedge 1\right)|w|^{-d-sp}\,dw\,dx\\
		&\leq C^{sp}\int_{Q_{n}^{*}}|f(x)|^{p}\int_{\R^{d}}\left(|z|^{p}\wedge 1\right)|z|^{-d-sp}l(Q_{n})^{-sp}\,dz\,dx\\
		&=C'\|f\|^{p}_{L^{p}(Q_{n}^{*})}l(Q_{n})^{-sp},
		\end{align*}
		with $C'$ depending on $s$, $d$, $p$ only.
	\end{proof}
	
	\begin{thm}\label{thm:finite}
		Suppose that $\Omega\subset\R^{d}$, $f\in W^{s,p}(\Omega)$ and
		\begin{equation}
		\int_\Omega \frac{|f(x)|^p}{\gamma(x)^{sp}}\,dx < \infty,
		\end{equation}
		where $\gamma(x)=\dist(x,\Omega^c)$. Then
		\begin{displaymath}
		\lim_{k\to \infty}[P^{\ep_k}f-f]_{W^{s,p}(\Omega)}=0,
		\end{displaymath}
		provided $\lim_{k\to\infty} \ep_k(Q)=0$ for every $Q\in\W$.
	\end{thm}
	
	\begin{proof}
		By \eqref{eq:Pandgn} and $\displaystyle\lim_{k\to\infty}\|\tau_{\ep_k(Q_{n})u}\left(g_{n}\right)-g_{n}\|^{p}_{L^{p}(\R^{2d})}=0$,
		it is enough to justify applications of Lebesgue dominated convergence theorem in Lemma~\ref{lem:gn}.
		To this end, we observe that
		\[
		\int_{\R^{d}}\|\tau_{\ep_k(Q_{n})u}\left(g_{n}\right)-g_{n}\|^{p}_{L^{p}(\R^{2d})}h(u)\,du\leq 2^{p}\|g_{n}\|^{p}_{L^{p}(\R^{2d})}.
		\] 
		Furthermore,
		\[
		\sum_{n=1}^{\infty} [f]^{p}_{W^{s,p}(Q_{n}^{*})} \leq 12^{d}[f]_{W^{s,p}(\Omega)}^{p}<\infty
		\]
	        and, by Whitney decomposition properties,
                $l(Q_{n})\geq\frac{\gamma(x)}{(5+\varepsilon)\sqrt{d}} \geq \frac{\gamma(x)}{6\sqrt{d}}$ for $x\in Q_{n}^{*}$, thus,
	\[
	\sum_{n=1}^{\infty}\|f\|^{p}_{L^{p}(Q_{n}^{*})} l(Q_{n})^{-sp} \leq\left(6\sqrt{d}\right)^{sp}\sum_{n=1}^{\infty}\int_{Q_{n}^{*}}\frac{|f(x)|^{p}}{\gamma(x)^{sp}}\,dx \leq\left(6\sqrt{d}\right)^{sp} 12^{d}\int_{\Omega}\frac{|f(x)|^{p}}{\gamma(x)^{sp}}\,dx<\infty.\qedhere
	\]
	\end{proof}
	
	We recall a geometric notion from  \cite{MR927080}.
	\begin{defn}
		A set $A\subset \R^d$ is {\em $\kappa$-plump}
		with $\kappa\in (0,1)$ if, for each $0<r< \diam(A)$ and each $x\in \bar{A}$, there
		is $z\in \bar B(x,r)$ such that
		$B(z,\kappa r)\subset A$.\end{defn}
	
	\begin{cor}
		Suppose that $\Omega\subset \R^d$ is an open set such that its complement $\Omega^c$ is $\kappa$-plump with some $\kappa\in (0,1)$, and $|\partial \Omega|=0$. Let $f\in W^{s,p}(\R^d)$ with $f=0$ on $\Omega^c$. Then
		\begin{displaymath}
		\lim_{k\to \infty}[P^{\ep_k}f-f]_{W^{s,p}(\R^d)}=0,
		\end{displaymath}
		provided $\lim_{k\to\infty} \ep_k(Q)=0$ for every $Q\in\W$.
	\end{cor}
	\begin{proof}
		We will show that such a~function $f$ satisfies the assumptions of Theorem~\ref{thm:finite} with the set $\R^d\setminus \partial\Omega$ in place of $\Omega$.
		Indeed, thanks to our assumptions we have $\int_{\Omega}\int_{\R^d\setminus\Omega}|f(x)|^{p}|x-y|^{-d-sp}\,dy\,dx<\infty$.
                Fix $R<\diam(\Omega^c)$ and let $x\in \Omega$ with $\gamma(x)=\dist(x,\Omega^c)<R$. Then
                \begin{align*}
		  \int_{\R^d\setminus\Omega}\frac{dy}{|x-y|^{d+sp}}&\geq
                  \int_{B(x,2\gamma(x)) \cap \Omega^c} \frac{dy}{|x-y|^{d+sp}}\\
                  &\geq C \gamma(x)^{-d-sp} |B(x,2\gamma(x)) \cap \Omega^c| \geq C' \gamma(x)^{-sp},
		\end{align*}
                where the last inequality follows from the $\kappa$-plumpness of $\Omega^c$.
                Thus
                \[
                \int_{ \{x \in \Omega:\,\gamma(x)<R\}} |f(x)|^p \gamma(x)^{-sp}\,dx < \infty.
                \]
                  Since $f\in L^p(\R^d)$ and $f=0$ on $\Omega^c$, it follows that  $\int_{\R^d\setminus\partial\Omega}\frac{|f(x)|^{p}}{\gamma(x)^{sp}}\,dx<\infty.$
	\end{proof}

        The next result is essentially a~fractional counterpart of the Meyers--Serrin theorem. The proof may be found for example in \cite[Theorem 3.25]{McLean}. We nevertheless provide the proof using our notation, as it is going to be modified in the next section.
        
	\begin{thm}[\cite{McLean}]\label{thm:Wsp}
		Let $p\in[1,\infty)$ and $s\in(0,1)$. Then the functions of a class $\mathcal{C}^{\infty}(\Omega) \cap W^{s,p}(\Omega)$ are dense in $W^{s,p}(\Omega)$.
	\end{thm}	
	\begin{proof}
		Let us fix a~function $f\in W^{s,p}(\Omega)$. Using notation \eqref{def:gn} from Lemma~\ref{lem:gn}, for all natural numbers $k$ and $n$, we choose $\ep_k(Q_n) < \frac{\varepsilon}{2k} l(Q_n)$ small enough so that the following inequality holds,
		\[
		\|\tau_t\left(g_{n}\right)-g_{n}\|^{p}_{L^{p}(\R^{2d})} < \frac{1}{k\, 2^n}\, , \quad 0<t<\ep_k(Q_n).
		\]
		Then from Lemma~\ref{lem:gn} it follows that
		\begin{align*}
		[P^{\ep_k}f-f]^{p}_{W^{s,p}(\Omega)} &= M\sum_{n=1}^{\infty}\int_{\R^{d}}\|\tau_{\ep_k(Q_{n})u}\left(g_{n}\right)-g_{n}\|^{p}_{L^{p}(\R^{2d})}h(u)\,du \\
		&\leq M \sum_{n=1}^{\infty} \frac{1}{k\, 2^n} = \frac{M}{k} \to 0,
		\end{align*}
		when $k\to \infty$. The convergence $P^{\ep_k}f \to f$ in $L^p(\Omega)$ follows from Theorem~\ref{thm:Lp},
		because $\ep_k(Q) \to 0$ for each $Q\in \W$.
		Finally, $P^{\ep_k}f \in C^\infty(\Omega)$ by Proposition~\ref{prop:cinfty}.
	\end{proof}	
	
	\section{Convergence in weighted spaces}\label{sec:weighted}
	In this section we extend our results to the case of weighted Sobolev spaces. Namely, for a~\emph{weight} $w$,
	i.e., a~nonnegative measurable function, we define the seminorm
	$$
	[f]_{W^{s,p}(\Omega,w)}=\left(\int_{\Omega}\int_{\Omega}\frac{|f(x)-f(y)|^{p}}{|x-y|^{d+sp}}w(y)w(x)\,dy\,dx\right)^{\frac{1}{p}},
	$$
	and the weighted $L^p$ norm
	\[
	\|f\|_{L^p(\Omega,w)} = \left(\int_\Omega |f(x)|^p w(x)\,dx \right)^{1/p}.
	\]
	
	We also denote
	$$
	\widetilde{W}^{s,p}(\Omega,w)=\left\{f\colon\Omega\to\mathbb{R}:\text{$f$ measurable, }\left[f\right]_{W^{s,p}(\Omega,w)}<\infty\right\}.
	$$
	
	\begin{prop}\label{prop:inclusion}
		If $w$ is locally comparable to a constant, that is for every compact $K\subset\Omega$ there is a constant $C_{K}>0$ such that $\frac{1}{C_{K}}\leq w(x)\leq C_{K}$ for all $x\in K$, then $\widetilde{W}^{s,p}(\Omega,w) \subset L^p_{loc}(\Omega)$.  
	\end{prop}
	\begin{proof}
		Fix two compact sets $K,L\subset\Omega$ of positive measure  and let $C=\sup_{x\in K}\sup_{y\in L}|x-y|<\infty$. To prove the inclusion, let us see that
		\begin{align*}
		\infty&>\int_{L}\int_{K}\frac{|f(x)-f(y)|^{p}}{|x-y|^{d+sp}}w(x)w(y)\,dx\,dy\geq C^{-d-sp}\int_{L}\int_{K}|f(x)-f(y)|^{p}w(x)w(y)\,dx\,dy.    
		\end{align*}
		By Fubini - Tonelli theorem, the inner integral $\int_{K}|f(x)-f(y)|^{p}w(x)\,dx$ is finite for almost all $f(y)$. Hence, for such $f(y)$, using the triangle inequality and the local boundedness of $w$, we have
		\[
		\int_{K}|f(x)|^{p}w(x)\,dx\leq 2^{p-1}\left(\int_{K}|f(x)-f(y)|^{p}w(x)\,dx+|f(y)|^{p}\int_{K}w(x)\,dx\right)<\infty.
		\]
		Now, 
		\[
		\int_{K}|f(x)|^{p}\,dx\leq C_{K}\int_{K}|f(x)|^{p}w(x)\,dx<\infty. \qedhere
		\]
	\end{proof}
	
	\begin{rem}\label{rem:continuous}
		If $w$ is continuous, then we can change $\Omega$ to $\Omega'=\Omega\setminus\{x:w(x)=0\}$. The set $\Omega'$ is still open and $w$ is locally comparable to a constant on $\Omega'$, so we can consider the space $W^{s,p}(\Omega',w)$ instead of $W^{s,p}(\Omega,w)$.
	\end{rem}
	\begin{lem}\label{lem:xy}
		If $y\in Q_{n}^{*}$ and $x\notin Q_{n}^{**}$, then $|x-y|\geq \frac{\varepsilon}{\varepsilon+\sqrt{d}}|x-x_{n}|$, where $x_{n}$ is the center of cube $Q_{n}$.
	\end{lem}
	\begin{proof}
		We have $|x-y|\geq \frac{\left(1+\varepsilon\right)^{2}l(Q_{n})-\left(1+\varepsilon\right)l(Q_{n})}{2}=\frac{\varepsilon(1+\varepsilon)l(Q_{n})}{2}$ and $|y-x_{n}|\leq\diam Q_{n}^{*}/2=(1+\varepsilon)l(Q_{n})\sqrt{d}/2$. Hence, $|x-y|\geq\frac{\varepsilon(1+\varepsilon)}{2}\frac{2}{(1+\varepsilon)\sqrt{d}}|y-x_{n}|=\frac{\varepsilon}{\sqrt{d}}|y-x_{n}|$.
		The assertion of the lemma follows from triangle inequality $|x-x_n| \leq |x-y| + |y-x_n|$.
	\end{proof}
	\begin{thm}\label{thm:weighted-dense} Suppose that $w$ is locally comparable to a constant or continuous and satisfies the condition
		\begin{equation}
		\label{conditon}
		\int_{\Omega}\frac{w(x)}{(1+|x|)^{d+sp}}\,dx<\infty
		\end{equation}
		Then $C^\infty(\Omega) \cap \widetilde{W}^{s,p}(\Omega,w)$ is dense in $\widetilde{W}^{s,p}(\Omega,w)$.
	\end{thm}
	\begin{proof}
		We extend $w$ to be $0$ outside $\Omega$. If $w$ is continuous, then we use Remark \ref{rem:continuous} and change $\Omega$ to $\Omega'$ in all the computations below. Similarly as in the previous cases, using the notations \eqref{def:gn} from Lemma \ref{lem:gn}, we have

		\begin{align*}
		[P^{\eta_{k}}f-f]^{p}_{W^{s,p}(\Omega,w)}\leq M\sum_{n=1}^{\infty}\int_{\mathbb{R}^{d}}\|\tau_{\ep_k(Q_{n})u}(g_{n})-g_{n}\|^{p}_{L^{p}(\mathbb{R}^{2d},w\times w)}h(u)\,du.
		\end{align*}
		We obtain for $t<\eta_k(Q_{n})$,
		\begin{align*}
		&\|\tau_t(g_{n})-g_{n}\|^{p}_{L^{p}(\mathbb{R}^{2d},w\times w)}\\
		&\leq\int_{Q_{n}^{*}}\int_{Q_{n}^{**}}\frac{|f(x-t)\psi_{n}(x-t)-f(y-t)\psi_{n}(y-t)-f(x)\psi_{n}(x)+f(y)\psi_{n}(y)|^{p}}{|x-y|^{d+sp}}w(y)w(x)\,dy\,dx\\
		&+2\int_{Q_{n}^{*}}\int_{\Omega\setminus Q_{n}^{**}}\frac{|f(x)\psi_{n}(x)-f(x-t)\psi_{n}(x-t)|^{p}}{|x-y|^{d+sp}}w(y)w(x)\,dy\,dx\\
		&=:I_{1}+2I_{2}.
		\end{align*}
		For the integral $I_{1}$ we have the following estimate
		\begin{align*}
		I_{1}&\leq C_{n}^2\int_{Q_{n}^{*}}\int_{Q_{n}^{**}}\frac{|f(x-t)\psi_{n}(x-t)-f(y-t)\psi_{n}(y-t)-f(x)\psi_{n}(x)+f(y)\psi_{n}(y)|^{p}}{|x-y|^{d+sp}}\,dy\,dx\\
		&\leq C_{n}^2\|\tau_t(g_{n})-g_{n}\|^{p}_{L^{p}(\mathbb{R}^{2d})},
		\end{align*}
		where $C_{n}=\displaystyle\sup_{x\in Q_{n}^{**}}w(x).$
		Let us now focus on the integral $I_{2}$. Using Lemma \ref{lem:xy}, if $x\in Q_n^{*}$ and $y\notin Q_n^{**}$, then $|x-y|\geq c|y-x_{n}|$ for $c=\varepsilon/(\varepsilon+\sqrt{d})$, when $x_{n}$ is the center of the cube $Q_n$. Thus, we obtain
		\begin{align*}
		I_{2}&\leq c^{-d-sp}\int_{Q_{n}^{*}}\int_{\Omega\setminus Q_{n}^{**}}\frac{|f(x)\psi_{n}(x)-f(x-t)\psi_{n}(x-t)|^{p}}{|y-x_{n}|^{d+sp}}w(y)w(x)\,dy\,dx\\
		&\leq C_{n}c^{-d-sp}\int_{Q_{n}^{*}}|f(x)\psi_{n}(x)-f(x-t)\psi_{n}(x-t)|^{p}\,dx \int_{\Omega\setminus Q_{n}^{**}}\frac{w(y)}{|y-x_n|^{d+sp}}\,dy\\
		&\leq D_{n}\|\tau_t\left(f\psi_{n}\right)-f\psi_{n}\|^{p}_{L^{p}(\R^{d})},
		\end{align*}
		where, thanks to Proposition \ref{prop:inclusion} the norm above is finite and
		\[
		D_{n}=\displaystyle C_n c^{-d-sp} \int_{\Omega\setminus Q_{n}^{**}}\frac{w(y)}{|y-x_n|^{d+sp}}\,dy.
		\]
		The integral above is finite, because for $y\not\in Q_n^{**}$ it holds
		$|y-x_n| \geq l(Q_n)/2$ and $|y-x_n| \geq |y| - |x_n|$, therefore $|y-x_n|$ is bounded from below by a~constant multiple of $1+|y|$.

		Now we need to repeat the proof of Theorem \ref{thm:Wsp}: for all natural numbers $k$ and $n$, we choose $\eta_{k}(Q_{n})<\frac{\varepsilon}{2k}l(Q_{n})$ such that 
		
		$$
		\|\tau_t\left(g_{n}\right)-g_{n}\|^{p}_{L^{p}(\R^{2d})}<\frac{1}{k2^{n+1}C_{n}^{2}}
		$$
		and
		$$
		\|\tau_t\left(f\psi_{n}\right)-f\psi_{n}\|^{p}_{L^{p}(\R^{d})}<\frac{1}{k2^{n+2} D_{n}},
		$$
		for $0<t<\eta_{k}(Q_{n}).$
		Hence, 
		\begin{align*}
		[P^{\eta_{k}}f-f]^{p}_{W^{s,p}(\Omega,w)}&\leq M\sum_{n=1}^{\infty}\int_{\mathbb{R}^{d}}\|\tau_{\ep_k(Q_{n})u}(g_{n})-g_{n}\|^{p}_{L^{p}(\mathbb{R}^{2d},w\times w)}h(u)\,du\\
		&\leq\frac{M}{k}\rightarrow 0, 
		\end{align*}
		
		when $k\rightarrow\infty.$

	\end{proof}
	\begin{thm}\label{thm:dense}
		Suppose that $w$ is locally comparable to a constant or continuous. Then
		$C^\infty(\Omega) \cap  L^{p}(\Omega,w)$ is dense in $L^{p}(\Omega,w)$.
	\end{thm}
	\begin{proof}
		If $w$ is continuous, then, according to Remark \ref{rem:continuous} we should replace $\Omega$ by $\Omega'$. Analogously as in the proof of Theorem~\ref{thm:Lp} we obtain that
		$$
		\|P^{\eta_{k}}f-f\|_{L^{p}(\Omega,w)}^{p}\leq M\sum_{n=1}^{\infty}C_{n}\int_{\mathbb{R}^{d}}\|\tau_{\ep_k(Q_{n})u}(f\psi_{n})-f\psi_{n}\|^{p}_{L^{p}(\Omega,w)}h(u)\,du.
		$$
		Since the function $\tau_{\ep_k(Q_{n})u}(f\psi_{n})$ has support in $Q_n^{**}$
		for $u\in \supp h$, taking $C_n=\displaystyle\sup_{x\in Q_{n}^{**}} w(x)$
		we obtain
		\[
		\|\tau_{\ep_k(Q_{n})u}(f\psi_{n})-f\psi_{n}\|^{p}_{L^{p}(\Omega,w)} \leq C_n \|\tau_{\ep_k(Q_{n})u}(f\psi_{n})-f\psi_{n}\|^{p}_{L^{p}(\R^d)}.
		\]
		We proceed as in the proof of Theorem~\ref{thm:weighted-dense} by
		choosing $\eta_{k}(Q_{n})<\frac{\varepsilon}{2k}l(Q_{n})$ such that $\|\tau_t(f\psi_{n})-f\psi_{n}\|^{p}_{L^{p}(\mathbb{R}^{d})}<\frac{1}{k2^{n+1}C_{n}}$ for $0<t<\eta_{k}(Q_{n})$ and we obtain the desired result.
	\end{proof}
	\begin{rem}
		
		Suppose that $\Omega=\R^{d}\setminus\{0\}$ and $w(x)=|x|^{-a}.$ The condition \eqref{conditon} becomes
		\begin{align*}
		\int_{\R^{d}\setminus\{0\}}\frac{dx}{|x|^{a}\left(1+|x|\right)^{d+sp}}& < \infty,
		\end{align*}
		which is equivalent to
		$$
		a\in(-sp,d).
		$$
		Analogous, but slightly different weighted Sobolev spaces were considered in  \cite{MR3420496}. Dipierro and Valdinoci considered density of compactly supported smooth functions in weighted Sobolev space $\dot{W}^{s,p}(\R^{d})=\widetilde{W}^{s,p}(\R^{d},w)\cap L^{p^{*}_{s}}(\R^{d},|\cdot|^{-2a/p})$ for $a\in [0,\frac{d-sp}{2})$ and $p_{s}^{*}=\frac{dp}{d-sp}$. Notice that however we do not have density of compactly supported functions, Theorems \ref{thm:weighted-dense} and \ref{thm:dense} combined provide a larger scale of the parameter $a$ and a general exponent $q$ instead of $p_{s}^{*}$. We can also change $\R^{d}\setminus\{0\}$ for any open set $\Omega$.
	\end{rem} 
	
	\section{Appendix}
	In this section we show how to generalise the results to the case of Sobolev spaces defined by some kernel $K$, see below.
	
	\begin{thm}\label{thm:K}
		Let $p\in[1,\infty)$, $\Omega\subset\R^{d}$ be an open set and let $K\colon[0,\infty)\to[0,\infty)$ be a measurable function such that 
		\[
		\int_{0}^{\infty}\left(x^p\wedge1\right)K(x)x^{d-1}\,dx<\infty.
		\]
		Denote
		\[
		[f]_{K} := \left( \int_\Omega \int_\Omega |f(x)-f(y)|^p K(|x-y|)\,dy\,dx \right)^{1/p}
		\]
		and consider the space
		\[
		X(\Omega)=\left\{f\in L^{p}(\Omega): [f]_K <\infty\right\}
		\]
		with the norm
		\[
		\|f\|_{X(\Omega)}=\left(\int_{\Omega}|f(x)|^{p}\,dx+ [f]_K^p \right)^{\frac{1}{p}}.
		\]
		Then the functions of a class $\mathcal{C}^{\infty}(\Omega) \cap X(\Omega)$ are dense in $(X(\Omega),\|\cdot\|_{X(\Omega)})$.
	\end{thm}
	
	\begin{proof}
		First we go through the proof of Lemma~\ref{lem:gn}, where we now estimate the seminorm $[P^{\ep_k}f-f]^{p}_K$
		and take
		\[
		g_{n}(x,y)=\left(f(x)\psi_{n}(x)-f(y)\psi_{n}(y)\right)K(|x-y|)^{\frac{1}{p}}, \quad \text{for $x,y\in\Omega$.}
		\]
		The only part of the proof that essentially changes is the estimate of $I_1$, which becomes
		\begin{align*}
		I_1&=\int_{Q_{n}^{*}}\int_{\Omega-x} |f(x)|^{p}|\psi_{n}(x)-\psi_{n}(x+w)|^{p} K(|w|) \,dw\,dx\\
		&=C' \|f\|^{p}_{L^{p}(Q_{n}^{*})} \int_{\R^{d}}\left(\frac{C^{p}|w|^{p}}{l(Q_{n})^{p}}\wedge 1\right) K(|w|)\,dw\,dx.
		\end{align*}
		We observe that
		\[
		\int_{\R^{d}}\left(\frac{C^{p}|w|^{p}}{l(Q_{n})^{p}}\wedge 1\right) K(|w|)\,dw \leq
		\left(\frac{C^{p}}{l(Q_{n})^{p}}\vee 1\right) \int_{\R^{d}}\left(|w|^{p} \wedge 1\right) K(|w|)\,dw < \infty.
		\]
		Having established an analogous version of Lemma~\ref{lem:gn}, we proceed as in the proof of Theorem~\ref{thm:Wsp} and obtain the desired result.
	\end{proof}


\begin{thebibliography}{1}
		
		\bibitem{MR3989177}
		A.~Baalal and M.~Berghout.
		\newblock Density properties for fractional {S}obolev spaces with variable
		exponents.
		\newblock {\em Ann. Funct. Anal.}, 10(3):308--324, 2019.
		
		\bibitem{MR3420496}
		S.~Dipierro and E.~Valdinoci.
		\newblock A density property for fractional weighted {S}obolev spaces.
		\newblock {\em Atti Accad. Naz. Lincei Rend. Lincei Mat. Appl.},
		26(4):397--422, 2015.
		
		\bibitem{MR3310082}
		A.~Fiscella, R.~Servadei, and E.~Valdinoci.
		\newblock Density properties for fractional {S}obolev spaces.
		\newblock {\em Ann. Acad. Sci. Fenn. Math.}, 40(1):235--253, 2015.
		
		\bibitem{MR3605166}
		H.~Luiro and A.~V. V\"{a}h\"{a}kangas.
		\newblock Beyond local maximal operators.
		\newblock {\em Potential Anal.}, 46(2):201--226, 2017.

		\bibitem{McLean}
		W. McLean
		\newblock Strongly elliptic systems and boundary integral equations.
		\newblock Cambridge University Press, Cambridge, 2000. 
		
		\bibitem{MR0290095}
		E.~M. Stein.
		\newblock {\em Singular integrals and differentiability properties of
			functions}.
		\newblock Princeton Mathematical Series, No. 30. Princeton University Press,
		Princeton, N.J., 1970.
		
		\bibitem{MR927080}
		J.~V\"{a}is\"{a}l\"{a}.
		\newblock Uniform domains.
		\newblock {\em Tohoku Math. J. (2)}, 40(1):101--118, 1988.
		
		\bibitem{MR3280034}
		Y.~Zhou.
		\newblock Fractional {S}obolev extension and imbedding.
		\newblock {\em Trans. Amer. Math. Soc.}, 367(2):959--979, 2015.
		
	\end{thebibliography}
\end{document}